\documentclass[11pt,twoside,reqno]{amsart}  
\usepackage{amssymb,verbatim} 
\usepackage{color} 
\usepackage{caption} 
\usepackage{hyperref}  
\usepackage[english]{babel}
\usepackage{amsmath}
\usepackage{amsfonts}
\usepackage{blindtext}
\usepackage{tabularx}
\usepackage{tabulary}
\usepackage{multirow}
\usepackage{graphicx}
\usepackage[all]{xy}
\usepackage{tikz}
\usepackage{amsthm}
\usepackage{enumerate}
\usepackage{titlesec}
\usepackage{lipsum}
\usepackage{bbm}

\titleformat{\section}
  {\normalfont\scshape\filcenter}{\thesection}{1em}{}
  
\numberwithin{equation}{section}

\titleformat{\subsection}
{\normalfont\normalsize\bfseries}{\thesubsection.}{1em}{}

\theoremstyle{plain}
	\newtheorem{Thm}{Theorem}
	\newtheorem{Lem}[Thm]{Lemma}
	\newtheorem{Cor}[Thm]{Corollary}
	
	\newtheorem{Prop}[Thm]{Proposition}
	\newtheorem{Rem}[Thm]{Remark}
	\newtheorem*{theoremA}{\bfseries Theorem A}
	\newtheorem*{theoremB}{\bfseries Theorem B}
\begin{document}
\title[Hochschild cohomology ring]{The Hochschild cohomology ring of the numerical semigroup algebras of embedding dimension two}  
 
\author[N. Tran]{Nghia T. H. Tran}  
\address{School of Mathematics, Statistics and Applied Mathematics, National University of Ireland, Galway, Ireland}  
\email{n.tranthihieu1@nuigalway.ie, emil.skoldberg@nuigalway.ie}  
\urladdr{}  

\author[E. Sk\" oldberg]{Emil Sk\" oldberg}  
\address{}  
\email{}  
\urladdr{} 
  
\thanks{Version: \today}  

\maketitle
\begin{abstract}
Let $a$ and $b$ be two coprime positive integers and $k$ an arbitrary field. We determine the ring structure of the Hochschild cohomology of the numerical semigroup algebras  $k[s^{a},s^{b}]$ of embedding dimension two (thus also complete intersections) in terms of generators and relations. In addition, we compute the Hilbert series for this cohomology ring.
\end{abstract}
\section{Introduction}
The ring structure of the Hochschild cohomology of an algebra is not known in general. Also, there are not many non-trivial examples of how to calculate this structure because of the complexity in describing the multiplicative structure. In \cite{Hol00}, Holm provided a description in terms of generators and relations of the Hochschild cohomology ring of the $k$-algebra $k[X]/\langle f \rangle$ where $f$ is a monic element of the polynomial ring $k[X]$ in a single variable. Of particular interest to us are cases of algebras with more variables; however, the method in \cite{Hol00} can only work in case of algebras with one variable. In \cite{TS18}, the authors provided a concrete approach to deal with computing the Hochschild cohomology ring of square-free monomial complete intersections in $n$ variables. The aim of this article is to obtain a description of the Hochschild cohomology ring of the numerical semigroup algebras $k[s^{a},s^{b}]\subseteq k[s]$ of embedding dimension two, a case of non-monomial complete intersections in more than one variable.

Our approach relies on the construction of the free resolution of complete intersections given by Guccione et al. \cite{GGRV}. We then provide a description of the Hochschild cohomology as a $k$-module by splitting the cochain complex into sub-complexes based on the features of cocycles. For the multiplicative structure, we interpret the cup product in terms of the Yoneda product. In order to compute the formula of the cup product of two elements in the module, starting from a cocycle we construct a chain map between the shifted resolution and the resolution itself. Building on previous work by one of the authors \cite{ES} on algebraic discrete Morse theory, we work out in Proposition \ref{contracting} an explicit description of the contracting homotopy, which allows us to construct the lifting map by combining the differentials and the contracting homotopy.

The formula of the differentials depends significantly on the relation of the two numbers $a$, $b$ and the characteristic $\mathrm{char}(k)$ of the field $k$. This yields that the structure of the Hochschild cohomology of $k[s^a,s^b]$ does the same. Therefore, we will consider this structure in two separate cases. The first case is that neither $a$ nor $b$ is divisible by $\mathrm{char}(k)$; hence the second is that $\mathrm{char}(k)$ is a divisor of $a$ or $b$, where we assume without loss of generality that $\mathrm{char}(k)$ is a divisor of $a$. For each of these two cases, we provide a description in terms of generators and relations of the Hochschild cohomology of $k[s^a,s^b]$ in the following theorems, which we shall prove in Subsections \ref{proofA} and \ref{proofB}. Subsequently, we calculate the Hilbert series of the Hochschild cohomology ring.

\begin{theoremA}[$\mathrm{char}(k) \nmid a,b$]\label{main.theorem}
Let $k$ be a field and $k[s^a,s^b]$ the numerical semigroup algebra, where $a$ and $b$ are coprime positive integers. The Hochschild cohomology algebra of $k[s^a,s^b]$ is isomorphic to the quotient ring
$$k[X_1,X_2,Y_1,Y_2,T]/\mathcal{I},$$
where $k[X_1,X_2,Y_1,Y_2,T]$ is a weighted graded commutative polynomial ring in which $\mathrm{deg}(X_1)=\mathrm{deg}(X_2)=0$, $\mathrm{deg}(Y_1)=\mathrm{deg}(Y_2)=1$ and $\mathrm{deg}(T)=2$; $\mathrm{wt}(X_1)=a$, $\mathrm{wt}(X_2)=b$, $\mathrm{wt}(Y_1)=0$, $\mathrm{wt}(Y_2)=ab-a-b$ and $\mathrm{wt}(T)=-ab$; and the ideal $\mathcal{I}$ is generated by the following relations: $X_1^b-X_2^a$, $X_1^{b-1}T$, $X_2^{a-1}T$, $Y_2T$, $Y_1^2$, $Y_2^2$, $Y_1Y_2$, $X_1Y_2-X_2^{a-1}Y_1$, $X_2Y_2-X_1^{b-1}Y_1$.
\end{theoremA}
\begin{theoremB}[$\mathrm{char}(k) \mid a$]\label{main.theorem2}
Let $k$ be a field with characteristic $n$ and $a$, $b$ two coprime integers such that $n$ is a divisor of $a$. Then the Hochschild cohomology algebra of $k[s^a,s^b]$ is isomorphic to the quotient ring
$$k[X_1,X_2,Y,T]/\mathcal{I},$$
where $k[X_1,X_2,Y,T]$ is a weighted graded commutative polynomial ring in which $\mathrm{deg}(X_1)=\mathrm{deg}(X_2)=0$, $\mathrm{deg}(Y)=1$ and $\mathrm{deg}(T)=2$; $\mathrm{wt}(X_1)=a$, $\mathrm{wt}(X_2)=b$, $\mathrm{wt}(Y)=-b$ and $\mathrm{wt}(T)=-ab$; and the ideal $\mathcal{I}$ is generated by the relations: 
\begin{itemize}
\item $X_1^{b}-X_2^{a}$, $X_1^{b-1}T$, and $Y^2-X_2^{a-2}T$ if $\mathrm{char}(k)=2$ and $4 \nmid a$; or
\item $X_1^{b}-X_2^{a}$, $X_1^{b-1}T$, and $Y^2$ otherwise.
\end{itemize}
\end{theoremB}
\textbf{Organization of the article.} In Section \ref{s1}, we introduce some auxiliary results about the semigroup and its Frobenius number, which are necessary for the later results.  In Section \ref{s2}, we review some basic facts on the topic and we describe the Hochschild cohomology module. After that, we divide the rest of the article into two independent parts according to the characteristic $\mathrm{char}(k)$ of the field $k$. Thus, in Section \ref{s3}, we treat the case in which neither $a$ nor $b$ is divisible by $\mathrm{char}(k)$, and in Section \ref{s4}, we treat the case that $\mathrm{char}(k)$ is a divisor of $a$.  In each of these two sections, we give a formula for the cup product and a proof of Theorems A and B above. Finally, we conclude the article by computing the Hilbert series of the algebras, which can be found in Subsections \ref{hilbert1} and \ref{hilbert2}.
 this article.  In Section \ref{s2}, we review some basic facts on the topic and especially, we calculate the Hochschild cohomology module, denoted by $\mathrm{HH}^*(A)$, based on the resolution introduced by Guccione et al. in \cite{GGRV}. After this section, we divide the rest of the article into two independent parts according to the characteristic $\mathrm{char}(k)$ of the field $k$. Thus, in Section \ref{s3}, we treat the case in which neither $a$ nor $b$ is divisible by $\mathrm{char}(k)$ and Section \ref{s4} for the case that $\mathrm{char}(k)$ is a divisor of $a$ (similarly for $b$).  In each of these two sections, we will make the chain map explicit based on the work of E. Sk\"oldberg \cite{ES} about algebraic discrete Morse theory. We then describe the multiplication on $\mathrm{HH}^*(A)$ which gives $\mathrm{HH}^*(A)$ a structure of a $k$-algebra. An explicit formula of multiplication is given in Corollary \ref{cup_product} and \ref{cup_product2}. Then we describe this algebra structure in term of generators and relations (Theorem \ref{main.theorem} and \ref{main.theorem2}) and we finish the article by providing the Hilbert series of the Hochschild cohomology algebra in subsections \ref{hilbert1} and \ref{hilbert2}.
\section{Some auxiliary results} \label{s1}
Let $S$ be the semigroup generated by $a$ and $b$, that is, $S:=\{ua+vb\mid u,v \in \mathbb{N}_0\}$. In this section, we prove some  results which will be used throughout the article.
\begin{Lem} \label{p2lemma1}
For an integer $d$, $db-a \in S$ if and only if $d \geq a$.
\proof ``$\Rightarrow$'': The Frobenius number of $S$ is $F(S) =ab-(a+b)$. Let $db-a\in S$ and suppose that $d=a-z$ where $z \in \mathbb{Z}, \ z \geq 1$. Thus $S\ni(a-z)b-a=F(S)-(z-1)b \Rightarrow F(S) \in S$ which is a contradiction. ``$\Leftarrow$'': $db-a = a(b-1)+b(d-a)\in S$ for any integer $d \geq a$.
\qed
\end{Lem}

To simplify the notation, we introduce $m_1=b(a-1)$ and $m_2=a(b-1)$ and define the sets 
$
S_i = \ \{\alpha \in \mathbb{Z} \ | \ \alpha-m_i \in S\}
$
for $i \in \{1,2\}$. Notice that $m_1, m_2 \in S$ (by Lemma \ref{p2lemma1}). The relevance of the next lemma will be seen later. 

\begin{Lem}\label{lemma2} The following are true:
\begin{enumerate}[(a)]
\item For $i \in \{1,2\}$, $S_i$ is equal to $\{m_i+\gamma \ | \ \gamma \in S\}$.
\item $S_1 \cap S_2$ is equal to $\{m_1+m_2\} \cup \{ab+\gamma \ | \ \gamma \in S\}$.
\end{enumerate}
\proof For (a), let $S'=\{m_i+\gamma \ | \ \gamma \in S\}$. We have $\alpha \in S_i \iff \alpha -m_i \in S \iff \exists \gamma \in S:\alpha-m_i=\gamma$, i.e., $\alpha=m_i+\gamma$. This is equivalent to $\alpha \in S'.$  For (b), let $S''=\{m_1+m_2\} \cup \{ab+\gamma \ | \ \gamma \in S\}$. Choose $\alpha \in S_1 \cap S_2$. In particular, $\alpha$ satisfies $\alpha-m_2 \in S$ and so we can write $\alpha = m_2+\beta$ for some $\beta = ua+vb \in S$ where $u,v \in \mathbb{N}_0$.
\begin{itemize}
\item If $u=0$ then $\alpha - {m_1} \in S \iff 
{-a+(v+1)b \in S}\iff v \geq a-1$ (by Lemma \ref{p2lemma1}). Thus $\alpha=m_1+m_2$ if $v=a-1$ and $\alpha=ab+\gamma$ where $ \gamma = vb-a\in S$ if $v \geq a$.
\item If $u>0$ then we can write $\beta = \gamma+a$ where $\gamma = (u-1)a+vb \in S$ giving $\alpha=ab+\gamma$.
\end{itemize}
Thus $S_1 \cap S_2 \subseteq S''$. The other inclusion is clear.
\qed
\end{Lem}
\section{A construction of Hochschild cohomology}\label{s2}
By setting $x_1\mapsto s^b$ and $x_2\mapsto s^a$, we have an isomorphism between algebras, $k[s^a,s^b]\cong \dfrac{k[x_1,x_2]}{\langle x_1^{a}-x_2^{b}\rangle}$. We will use both algebras according to our convenience. We now interpret the minimal resolution given by Guccione et al. (see \cite{GGRV}) for the case of the quotient ring of $k[x_1,x_2]$ modulo $\langle x_1^{a}-x_2^{b}\rangle$, the ideal generated by the binomial $x_1^{a}-x_2^{b}$. Let us denote by $A^{op}$ the opposite algebra of $A$. The tensor product will be taken over $k$, i.e., $\otimes=\otimes_k$. We denote by $A^e:=A\otimes A^{op}$ the enveloping algebra of $A$. The following complex $\mathbf{F}$ is a free $A^e$-resolution of $A$:
	\begin{equation} \label{complex1}
\cdots \longrightarrow F_2 			\xrightarrow{\text{ $d_2$ }} F_1 \xrightarrow{\text{ $d_1$ }} F_0 \xrightarrow{\text{ $\epsilon$ }}A\longrightarrow 0,
	\end{equation}
where $F_m$ is the finitely generated free $A^e$-module with basis elements $e_{i_1\cdots i_r}\cdot t^{(q)}$ ($r,q \geq 0$ and $r+2q=m$), where by $e_{i_1\cdots i_r}$ or $e_I$ ($I=\{i_1,\ldots,i_r\}\subseteq \{1,2\}, i_1 < \cdots <i_r$), we mean $e_{i_1}\wedge \cdots \wedge e_{i_r}$. Then, we have that \eqref{complex1} is an exact sequence of free $A^e$-modules $F_m$ with
$$\epsilon:A^e\to A, \qquad a\otimes b \mapsto ab$$
and the differentials $d_m$ (briefly $d$) are defined as follows:
\begin{center}
$d(e_1)=s^b\otimes 1-1\otimes s^b$;\\
$d(e_2)=s^a\otimes 1-1\otimes s^a$;\\
$d(t)=\sum\limits_{i=0}^{a-1}{s^{ib}\otimes s^{b(a-1-i)}\cdot e_1}-\sum\limits_{i=0}^{b-1}{s^{ia}\otimes s^{a(b-1-i)}\cdot e_2};$
\end{center}
inductively,
\begin{center}
$d(xy)=d(x)y+(-1)^{\alpha}xd(y)$ where $x \in F_\alpha$.
\end{center}
Alternatively, we can write 
\begin{equation*} d(e_It^{(q)})=
\left\{
	\begin{array}{ll}
		d(t)t^{(q-1)},  & I = \emptyset, \\
		d(e_1)t^{(q)}-e_1d(t)t^{(q-1)},  & I = \{1\}, \\
		d(e_2)t^{(q)}-e_2d(t)t^{(q-1)},  & I = \{2\}, \\
		\left(d(e_1)e_2-d(e_2)e_1\right)t^{(q)},  & I = \{1,2\}.
	\end{array}
\right. (\text{if } q>0)
\end{equation*}
Applying the contravariant functor $\mathrm{Hom}_{A^e}(-,A)$ to the truncation of the above resolution, we obtain a new complex: 
\begin{equation*}
0 \longrightarrow \text{Hom}_{A^e}(F_0,A)  \xrightarrow{\text{ }\text{ }\text{ }d^1 \text{ }\text{ }\text{ }} \text{Hom}_{A^e}(F_1,A) \xrightarrow{\text{ }\text{ }\text{ }d^2\text{ }\text{ }\text{ }} \text{Hom}_{A^e}(F_2,A) \longrightarrow\cdots.
\end{equation*}
The $i$-th Hochschild cohomology of $A$ is the module
$$\mathrm{HH}^i(A):=\dfrac{\mathrm{Ker}(d^{i+1})}{\mathrm{Im}(d^i)},$$ where $d^0$ is taken to be the zero map. Now the Hochschild cohomology module of $A$ is defined to be the direct sum of these components, 
$
\mathrm{HH}^*(A):= \mathop \oplus \limits_{i\geq 0}\mathrm{HH}^i(A).
$

For $m\in \mathbb{N}_0$, let $\overline{F}_m$ be the free $k$-module generated by the same basis elements as $F_m$. Then, there is an isomorphism between the following $k$-spaces: 
\begin{equation*}
\text{Hom}_{A^e}(F_m,A) \cong  \text{Hom}_k(\overline{F}_m,A).
\end{equation*}
Thus one gets the following complex:
\begin{equation} \label{complex2}
0 \longrightarrow \text{Hom}_{k}(\overline{F}_0,A)  \xrightarrow{\text{ }\text{ }\text{ }\partial^1 \text{ }\text{ }\text{ }} \text{Hom}_{k}(\overline{F}_1,A) \xrightarrow{\text{ }\text{ }\text{ }\partial^2 \text{ }\text{ }\text{ }}\text{Hom}_{k}(\overline{F}_2,A)\longrightarrow \cdots,
\end{equation}
where  the differential $\partial$ will be described later.\\
Let $e_It^{(q)}$ be a basis element in $\overline F_m$ and $s^\alpha$ a basis element in $A$. Let $(e_It^{(q)},s^\alpha)$ be the $k$-linear map in $\mathrm{Hom}_k(\overline F_m,A)$ which sends $e_It^{(q)}$ to $s^\alpha$ and other basis elements to 0, that is,
\begin{equation*}
(e_It^{(q)},s^\alpha)(e_Jt^{(p)})=
\left\{
	\begin{array}{ll}
	s^\alpha & \text{if }J=I \text{ and } p=q,\\
	0 & \text{otherwise}.
	\end{array}
\right.
\end{equation*}
The set of all such $k$-linear maps is a $k$-basis of the module $\mathrm{Hom}_k(\overline F_m,A)$. We use the notation $[(e_It^{(q)},s^\alpha)]$ to denote the residue class represented by $(e_It^{(q)},s^\alpha)$ in $\mathrm{HH}^*(A)$. Now we are in the position to describe the formula for $\partial$.
\begin{Lem} \label{lemma.partial}
The homomorphism $\partial$ in \eqref{complex2} is given by:
\begin{equation*} \partial(e_It^{(q)},s^{\alpha})=
\left\{
	\begin{array}{ll}
		0  & \text{if } I = \emptyset, \\
		a(t^{(q+1)}, s^{\alpha+m_1})  & \text{if } I = \{1\}, \\
		-b(t^{(q+1)}, s^{\alpha+m_2}) &  \text{if } I = \{2\}, \\
		b(e_1t^{(q+1)}, s^{\alpha+m_2})+a(e_2t^{(q+1)}, s^{\alpha+m_1})  &  \text{if } I = \{1,2\}.
	\end{array}
\right. 
\end{equation*}
\proof
From the following diagram
\begin{equation*}
\xymatrix{\mathrm{Hom}_{A^e}(F_m,A) \ar[d]_{\cong} \ar[r]^{d^{m+1}} & \mathrm{Hom}_{A^e}(F_{m+1},A)\ar[d]_{\cong}\\
\mathrm{Hom}_{k}(\overline{F}_m,A) \ar[u] \ar[r]^{\partial^{m+1}}&\mathrm{Hom}_{k}(\overline{F}_{m+1},A) \ar[u]}
\end{equation*}
we can derive $\partial^{m+1}$ from $d^{m+1}$ straightforwardly.\\
As $f=(e_It^{(q)},s^\alpha)$ is a basis element of $\mathrm{Hom}_{k}(\overline{F}_m,A)$, $f$ is identified with a function in $\mathrm{Hom}_{A^e}(F_m,A)$. A direct calculation shows that it is $(e_It^{(q)},s^\alpha)$, which is also denoted by $f$ by abuse of notation. We have the homomorphism
\begin{equation*}
\begin{array}{cccc}
d^{m+1}: & \mathrm{Hom}_{A^e}(F_m,A) & \longrightarrow & \mathrm{Hom}_{A^e}(F_{m+1}, A)\\
       & f                & \longmapsto &d^{m+1}(f):=f\circ d
\end{array}
\end{equation*}
%sends $f$ to $d^{m+1}(f)$, an $A^e$-homomorphism from $F_{m+1}$ to $A$.\\
For a basis element $e_Jt^{(r)} \in F_{m+1}$, we can compute $[d^{m+1}(f)](e_Jt^{(r)})$ directly and the result is summarized in the following table:
\begin{center}
\begin{tabulary}{\textwidth}{|C|C|}
\hline
{}& $[d^{m+1}(f)](e_Jt^{(r)})$\\
\hline
\multirow{3}{*} {$J = \emptyset$} & 
$as^{\alpha +m_1}$, if $I= \{1\}$ and $r-1=q$\\
{}& 
$-bs^{\alpha +m_2}$, if $I= \{2\}$ and $r-1=q$\\
{}& 0, otherwise\\\cline{2-2}
\hline
\multirow{2}{*}{$J = \{1\}$} & $bs^{\alpha+m_2}$, if $I= \{1,2\}$ and $r-1=q$\\
{}& 0, otherwise\\\cline{2-2}
\hline
\multirow{2}{*}{$J = \{2\}$} & $as^{\alpha+m_1}$, if $I= \{1,2\}$ and $r-1=q$\\
{}& 0, otherwise\\\cline{2-2}
\hline
$J = \{1,2\}$ & 0\\
\hline 
\end{tabulary}
\end{center}
In other words, we have the formula of $\partial$ as desired.
\qed
\end{Lem}
We will now consider two cases. The first, Case I, is when the characteristic $\mathrm{char}(k)$ of the field $k$ is neither a divisor of $a$ nor of $b$, and the second, Case II, is when $\mathrm{char}(k)$ divides one of $a$ or $b$, which we without loss of generality assume to be $a$.
\section{The ring structure of $\mathrm{HH}^*(A)$ - Case I}\label{s3}
\subsection{The structure of $k$-module $\mathrm{HH}^*(A)$  via sub-complexes} \label{classification}
The $k$-vector space $\mathrm{Hom}_k(\overline{F}_m,A)$ is generated by the basis elements of cohomological degree $m$:
\begin{equation*}
\mathrm{Hom}_k(\overline{F}_m,A)=\mathop\oplus\limits_{|I|+2q=m}{k(e_It^{(q)},s^\alpha)},
\end{equation*} 
where $k(e_It^{(q)},s^\alpha)$ is the $k$-vector space generated by $(e_It^{(q)},s^\alpha)$. To simplify the notation we will use the notation $k$ instead of $k(e_It^{(q)},s^\alpha)$, justified by the isomorphism $k(e_It^{(q)},s^\alpha)\cong k$. In order to describe the Hochschild cohomology, we split the cohomology complex into subcomplexes.
Let $\Gamma:=\{(t^{(q)},s^\alpha)\mid q \in \mathbb{N}_0,\alpha \in S\}$,  the set of all basis elements that is in the kernel of $\partial$. For each element $\gamma \in \Gamma$, we construct the complex $M_\gamma$ which includes $\gamma$ as the generator of the rightmost non-zero entry ($\cdots \longrightarrow \bullet \longrightarrow 0$). Each $M_\gamma$ is a subcomplex of \eqref{complex2}. Moreover, by Lemma \ref{lemma.partial} there are only four options for such $M_\gamma$ as showed in the following table:
\begin{center}
\begin{tabulary}{\textwidth}{|L|C|L|}
  \hline
  $q$ & Type of subcomplex & Condition(s) for $\alpha$ in $(t^{(q)},s^\alpha)$\\
  \hline 
  $q=0$ & $0 \longrightarrow k \longrightarrow 0$ &    \\
  \hline
  \multirow{3}{*}{$q=1$} & $0 \longrightarrow k \longrightarrow 0$ & $\left\{ {\begin{array}{*{20}{c}}
{\alpha - {m_1} \notin S}\\
{\alpha - {m_2} \notin S}
\end{array}} \right.$\\ \cline{2-3}
 & $0 \longrightarrow k \longrightarrow k \longrightarrow 0$ & $\alpha-m_1 \in S \text{ xor } \alpha-m_2 \in S$\\ \cline{2-3}
&$0 \longrightarrow k^2 \longrightarrow k \longrightarrow 0$ &$\left\{ {\begin{array}{*{20}{c}}
{\alpha - {m_1} \in S}\\
{\alpha - {m_2} \in S}
\end{array}} \right.$\\
\hline
\multirow{4}{*}{$q \geq 2$} & $0 \longrightarrow k \longrightarrow 0$ & $\left\{ {\begin{array}{*{20}{c}}
{\alpha - {m_1} \notin S}\\
{\alpha - {m_2} \notin S}
\end{array}} \right.$\\\cline{2-3}
 & $0 \longrightarrow k \longrightarrow k \longrightarrow 0$ & $\alpha-m_1 \in S \text{ xor } \alpha-m_2 \in S$ \\\cline{2-3}
&$0 \longrightarrow k^2 \longrightarrow k \longrightarrow 0$ &$\left\{ {\begin{array}{*{20}{c}}
{\alpha - {m_1} \in S}\\
{\alpha - {m_2} \in S}\\
{\alpha - {m_1} - {m_2} \notin S}
\end{array}} \right.$\\\cline{2-3}
&$0 \longrightarrow k \longrightarrow k^2 \longrightarrow k \longrightarrow 0$&$\alpha - {m_1} - {m_2} \in S$\\
\hline
\end{tabulary} 
\captionof{table}{Classification of the sub-complexes}
\end{center}
We can write in detail the above subcomplex as follows:
\begin{itemize}
\item [Type 1:] $0 \longrightarrow k(t^{(q)},s^\alpha)\longrightarrow 0$;
\item [Type 2:] $0 \longrightarrow k(e_1t^{(q-1)},s^{\alpha-m_1}) \longrightarrow k(t^{(q)},s^\alpha)\longrightarrow 0$; or \\
$0 \longrightarrow k(e_2t^{(q-1)},s^{\alpha-m_2}) \longrightarrow k(t^{(q)},s^\alpha)\longrightarrow 0$;
\item [Type 3:] $0 \longrightarrow k(e_1t^{(q-1)},s^{\alpha-m_1})\oplus k(e_2t^{(q-1)},s^{\alpha-m_2}) \longrightarrow k(t^{(q)},s^\alpha)\longrightarrow 0$; 
\item [Type 4:] $0 \longrightarrow k(e_1\wedge e_2t^{(q-2)},s^{\alpha-m_1-m_2}) \longrightarrow k(e_1t^{(q-1)},s^{\alpha-m_1})\oplus k(e_2t^{(q-1)},s^{\alpha-m_2}) \longrightarrow k(t^{(q)},s^\alpha)\longrightarrow 0$. 
\end{itemize}
\begin{Prop} \label{prop.cocycles}
The complex \eqref{complex2} can be written as the direct sum below:
\begin{equation*}
\mathrm{Hom}_k(\overline{F}_\bullet,A)=\mathop\oplus\limits_{\gamma \in \Gamma}{M_\gamma}.
\end{equation*}
\proof
By Lemma \ref{lemma.partial}, $M_\gamma$ is a subcomplex of \eqref{complex2}. Then we have the first inclusion $\mathop\oplus\limits_{\gamma \in \Gamma}{M_\gamma} \subseteq \mathrm{Hom}_k(\overline{F}_\bullet,A)$. Let us consider an arbitrary non-zero basis element $E=(e_It^{(q)},s^\alpha) \in \mathrm{Hom}_k(\overline{F}_m,A)$ for some $m \in \mathbb{N}_0$. If $I=\emptyset$, then $E \in \Gamma$ and $M_E$ is the subcomplex containing $E$. If $I \neq \emptyset$, we have the following cases:\\
(a) $I=\{1\}$:
\begin{itemize}
\item If $\alpha - m_2\in S$, then $E$ occurs in the subcomplex 
$$0 \longrightarrow k(e_1\wedge e_2t^{(q-1)},s^{\alpha-m_2}) \longrightarrow kE\oplus k(e_2t^{(q)},s^{\alpha +m_1-m_2}) \longrightarrow k(t^{(q+1)},s^{\alpha+m_1})\longrightarrow 0$$
\item If $\alpha - m_2\notin S$, then $E$ occurs in the subcomplex 
$$0 \longrightarrow kE \longrightarrow k(t^{(q+1)},s^{\alpha+m_1})\longrightarrow 0\text{, if } \alpha+m_1-m_2 \notin S;\text{ and}$$
$$0 \longrightarrow kE\oplus k(e_2t^{(q)},s^{\alpha+m_1-m_2}) \longrightarrow k(t^{(q+1)},s^{\alpha+m_1})\longrightarrow 0  \text{, if } \alpha+m_1-m_2 \in S.$$
\end{itemize}
(b) $I=\{2\}$, similarly.\\
(c) If $I=\{1,2\}$, then $E$ occurs in the subcomplex 
$$0 \longrightarrow kE \longrightarrow k(e_1t^{(q+1)},s^{\alpha+m_2})\oplus k(e_2t^{(q+1)},s^{\alpha+m_1}) \longrightarrow k(t^{(q+2)},s^{\alpha+m_1+m_2})\longrightarrow 0.$$
We see that any basis element $E$ is contained in a unique subcomplex which belongs to Type 1 to 4 as above. So the inverse inclusion is obtained and the result follows. 
\qed
\end{Prop}
\subsection{Classification of cocycles}
The following is a direct consequence of Lemma \ref{lemma.partial} and Proposition \ref{prop.cocycles}. 
\begin{Cor} \label{cor.cocycles}
The $k$-vector space $\mathop \oplus\limits_{i\in \mathbb{N}_0}\mathrm{Ker}(\partial^i)$ is generated by the following elements:
\begin{itemize}
\item $(t^{(q)},s^{\alpha})$ where $\alpha \in S$;
\item $b(e_1t^{(q)},s^{\alpha -m_1})+a(e_2t^{(q)},s^{\alpha -m_2})$ where $\alpha \in S$ such that $\alpha-m_1 \in S$ and $\alpha-m_2 \in S$.
\end{itemize}
\end{Cor}
We will call the elements described in the above corollary the \textit{standard} elements. In the following remarks, we will give more details about these elements.
\begin{Rem} \label{cocycles_details}
By Lemma \ref{lemma2}, for any $\alpha \in \mathbb{Z}$ such that $\alpha-m_1 \in S$ and $\alpha-m_2 \in S$ we have that
\begin{equation*}
\alpha \in \{m_1+m_2\} \cup \{ab+\gamma \ | \ \gamma \in S\}.
\end{equation*}
If $\alpha =m_1+m_2$, we have the cocycles $b(e_1t^{(q)},s^{m_2})+a(e_2t^{(q)},s^{m_1})$.\\
If $\alpha =ab+\gamma$ for some $\gamma \in S$, then we get the cocycles
$b(e_1t^{(q)},s^{\gamma +b})+a(e_2t^{(q)},s^{\gamma + a})$ where $\gamma \in S$, which is distinguished from the cocycles above.
\end{Rem}
\begin{Rem}\label{non_zero_cocycles} 
 According to Table 1 and Corollary \ref{cor.cocycles}, we are able to identify all standard cocycles that are the representatives of the non-zero cohomology classes in $\mathrm{HH}^*(A)$ as follows:
\begin{equation}
(1,s^{\alpha}), \text{ where } \alpha \in S
\end{equation}
\begin{equation} \label{cocycles_type2}
(t^{(q)},s^{\alpha}), \text{ where } \left\{ {\begin{array}{*{20}{c}}
{q>0}\\
{\alpha \in S}\\
{\alpha - {m_1} \notin S}\\
{\alpha - {m_2} \notin S}
\end{array}} \right.
\end{equation}
\begin{equation}\label{cocycles_3}
b(e_1t^{(q)},s^{\alpha -m_1})+a(e_2t^{(q)},s^{\alpha -m_2}), \text{ where }\left\{ {\begin{array}{*{20}{c}}
{\alpha \in S}\\
{\alpha - {m_1} \in S}\\
{\alpha - {m_2} \in S}\\
{\alpha - {m_1} - {m_2}\notin S \text{ if } q>0}
\end{array}} \right.
\end{equation}

We can express the cocycles in \eqref{cocycles_3} as all the elements of the set
\begin{equation*}
\{b(e_1t^{(q)},s^{\gamma+b})+a(e_2t^{(q)},s^{\gamma+a}) \ | \ \gamma \in S; \ \gamma-F(S) \not\in S \text{ if } q>0\}
\end{equation*}
together with the single cocycle $b(e_1,s^{m_2})+a(e_2,s^{m_1})$. Indeed, by Remark \ref{cocycles_details}, we have $\gamma+b=\alpha-m_1$. Then $\alpha-m_1-m_2 \notin S$ is equivalent to $\gamma+b-m_2=\gamma-(ab-a-b)=\gamma-F(S)\notin S$.
\end{Rem}
In the above remark, we have described a $k$-basis of $\mathrm{HH}^*(A)$. Next, we will construct a multiplicative structure on $\mathrm{HH}^*(A)$.
\subsection{A Morse matching on $\mathbf{F}$}
Let $\mathbf{F}$ be a free resolution of the $A^e$-module $A$ and $f:F_i \to A$ be an $A^e$-homomorphism such that $f\circ d_{i+1}=0$.  By the comparison theorem, there is a chain map $\tilde f$ consisting of homomorphisms $\tilde f_m$, $m \in \mathbb{N}_0$ that makes the following diagram commute, moreover such a chain map is unique up to chain homotopy.
\begin{equation} \label{diagram1}
 \xymatrix{
 \cdots \ar[r]& {F_{i+2}}\ar[d]^{\tilde f_2}  \ar[r]^{d_{i+2}} &F_{i+1}\ar[d]^{\tilde f_1}\ar[r]^{d_{i+1}} &F_i\ar[d]^{\tilde f_0} \ar[rd]^f\\
\cdots \ar[r]& {F_2}  \ar[r]^{d_2} &F_1\ar[r]^{d_1} &F_0\ar[r]^{\epsilon}&A}
 \end{equation}
 In theory, it is always possible to construct such a chain map. However, it is not easy in practice. The goal of this section is to provide an explicit chain map in case of our resolution that makes the above diagram commute. In more details, we will base ourselves on a work of Sk\"oldberg \cite{ES} to construct a contracting homotopy $\phi$ which consists of maps of degree 1 $\phi_{j}: F_j \to F_{j+1}$. The homomorphism $\tilde{f}$ is given by setting
\begin{equation*}
\tilde{f}_0(x)=1 \otimes f(x)
\end{equation*}
and for any $j>0$, we define $\tilde{f}_j$ inductively on the $A^e$-basis elements by
\begin{equation*}
\tilde{f}_j:=\phi_{j-1} \circ \tilde{f}_{j-1}\circ d_{i+j}
\end{equation*}
and extend linearly for other elements.\\
The chain map $\tilde{f}$ defined as above makes diagram \ref{diagram1} commute. In the next steps, we will make this chain map explicit.

To denote the elements in the algebra $k[x_1,x_2]$ and their cosets in the quotient ring $k[x_1,x_2]/\langle x_1^a-x_2^b \rangle$, we use the same notation if there are no ambiguities. Then the $k$-basis of $k[x_1,x_2]/\langle x_1^a-x_2^b \rangle$ consists of all elements $x_1^{u}x_2^{v}$, where $u\geq 0$ and $0 \leq v<b$. From now on, these are default conditions whenever we mention the elements in  $k[x_1,x_2]/\langle x_1^a-x_2^b \rangle$.

We can consider the complex $\mathbf{F}$ as a chain complex of $k\otimes A$-modules together with a direct sum decomposition as follows:
$$F_m=\mathop \oplus \limits _{\alpha \in \mathfrak{I}_m}F_{\alpha},$$
where $\{\mathfrak{I}_m\}_{m \in \mathbb{N}_0}$ is a family of mutually disjoint index sets given by
$$\mathfrak{I}_m=\{(u,v,I,q)\mid u\geq 0, 0 \leq v<b, 2|I|+q=m\}.$$
Here the index $(u,v,I,q)$ corresponds to the basis element $x_1^{u}x_2^{v}\otimes 1 \cdot e_It^{(q)}$ which generates the $k\otimes A$-module $F_{(u,v,I,q)}$. We write $d_{\beta,\alpha}$ for the component of $d$ going from $F_{\alpha}$ to $F_{\beta}$. Now $\mathbf{F}$ has the structure of a \textit{based complex}, see \cite{ES} Section 2. Let $G_{\mathbf{F}}$ be the digraph with vertex set $V=\mathop \cup \limits_{m \in \mathbb{N}_0}\mathfrak{J}_m$ and with a directed edge $\alpha\to \beta$ whenever the component $d_{\beta,\alpha}$  is non-zero. Next, we construct a partial matching $M$ on $\mathbf{F}$ by 
\begin{align*}
\left. {\begin{array}{*{20}{c}}
{x_1^{u}x_2^{v}\otimes 1 \cdot t^{(q)} \longrightarrow x_1^{u-1}x_2^{v}\otimes 1 \cdot e_1t^{(q)}}\\
{x_1^{u}x_2^{v}\otimes 1 \cdot e_2t^{(q)} \longrightarrow x_1^{u-1}x_2^{v}\otimes 1 \cdot e_1e_2t^{(q)}}
\end{array}} \right\}\text{ where }u>0, 0\leq v<b
\end{align*}
\begin{align*}
x_2^{v}\otimes 1 \cdot t^{(q)} &\longrightarrow x_2^{v-1}\otimes 1 \cdot e_2t^{(q)}, \text{ where } 0<v<b\\
x_2^{b-1}\otimes 1 \cdot e_2t^{(q)} &\longrightarrow 1\otimes 1 \cdot t^{(q+1)}.
\end{align*}
We denote by $G_{\mathbf{F}}^M$ the digraph with the same vertex set $V$ and the edge set obtained from $G_{\mathbf{F}}$ by reversing the direction of each arrow $\alpha \to \beta$ whenever  $\beta \to \alpha$ in $M$.
 For each edge $\alpha \to \beta$ in $M$, it is clear that the corresponding component of the differential $d_{\beta,\alpha}$ is an isomorphism. Now we only need to check that there are no directed cycles in $G_{\mathbf{F}}^M$ to see that $M$ is a Morse matching. By observing the formula of the differential $d$ and the matching $M$, we check the absence of directed cycles as follows:\\
(i) If we have a path
$$x_2^{b-1}\otimes 1 \cdot e_2t^{(q)} \longrightarrow 1\otimes 1 \cdot t^{(q+1)}\longrightarrow x_1^ux_2^v\otimes 1 \cdot e_It^{(r)}$$
in $G_{\mathbf{F}}^M$ where the two first vertices are matched, then one gets $I=\{1\}$ or ($I=\{2\}$, $u=0$ and $v<b-1$), i.e., this path ends here and hence, it cannot form a cycle.\\
(ii) Similarly, if we have a path
$$x_2^{v}\otimes 1 \cdot t^{(q)} \longrightarrow x_2^{v-1}\otimes 1 \cdot e_2t^{(q)}\longrightarrow x_1^mx_2^n\otimes 1 \cdot e_It^{(r)}(\text{where } 0<v<b),$$
then $I=\{1,2\}$ (i.e., the path must end here and there is no cycle formed) or one has $x_1^mx_2^n\otimes 1 \cdot e_It^{(r)}=x_2^{v-1}\otimes 1 \cdot t^{(q)}$. Thus, the path becomes
$$x_2^{v}\otimes 1 \cdot t^{(q)} \longrightarrow x_2^{v-1}\otimes 1 \cdot e_2t^{(q)}\longrightarrow x_2^{v-1}\otimes 1 \cdot t^{(q)}\longrightarrow x_2^{v-2}\otimes 1 \cdot e_2t^{(q)}\longrightarrow \cdots $$
where the power of $x_2$ is declining and the path eventually terminates at $1\otimes 1 \cdot t^{(q)}$. Thus, no cycle is formed by this path.\\
(iii) Let us consider the path 
$$x_1^{u}x_2^{v}\otimes 1 \cdot t^{(q)} \longrightarrow x_1^{u-1}x_2^{v}\otimes 1 \cdot e_1t^{(q)}\longrightarrow x_1^mx_2^n\otimes 1 \cdot e_It^{(r)} (\text{where } u>0).$$
Then we have either $I=\{1,2\}$ (i.e., the path ends here) or $x_1^mx_2^n\otimes 1 \cdot e_It^{(r)}=x_1^{u-1}x_2^{v}\otimes 1 \cdot t^{(q)}$ and we can extend this path as follows:
$$x_1^{u}x_2^{v}\otimes 1 \cdot t^{(q)} \longrightarrow x_1^{u-1}x_2^{v}\otimes 1 \cdot e_1t^{(q)}\longrightarrow x_1^{u-1}x_2^{v}\otimes 1 \cdot t^{(q)} \longrightarrow \cdots \longrightarrow x_2^{v}\otimes 1 \cdot t^{(q)}$$
and continue with the path in (ii), i.e., there is no directed cycle.\\
(iv) For the last one, the path
$$x_1^{u}x_2^{v}\otimes 1 \cdot e_2t^{(q)} \longrightarrow x_1^{u-1}x_2^{v}\otimes 1 \cdot e_1e_2t^{(q)} \longrightarrow x_1^mx_2^n\otimes 1 \cdot e_It^{(r)}(\text{where } u>0)$$
gives us either $I=\{1\}$ (which ends the path) or $x_1^mx_2^n\otimes 1 \cdot e_It^{(r)}=x_1^{u-1}x_2^{v}\otimes 1 \cdot e_2t^{(q)}$. By continuing this argument, this path is extended to
$$x_2^{v}\otimes 1 \cdot e_2t^{(q)},$$
which ends here if $v<b-1$ and ends at $1\otimes 1 \cdot t^{(q+1)}$ if $v=b-1$.
Hence, there is no directed cycle in $G_{\mathbf{F}}^M$ and $M$ is a Morse matching as desired.

We will now give the formula of the contracting homotopy $\phi$ (see \cite{ES} for the definition of $\phi$) for our case in the following proposition.
\begin{Prop} \label{contracting}
Let $x=x_1^ux_2^v\otimes 1 \cdot e_It^{(q)}$ be a basis element of the $k\otimes A$-complex $\mathbf{F}$. We then have the formula of $\phi$ as follows:\begin{itemize}
\item $I=\{1\}$ or $\{1,2\}$: $\phi(x)=0$;
\item $I=\{\emptyset\}$: $\phi(x)=\sum\limits_{i=0}^{u-1}{x_1^ix_2^v\otimes x_1^{u-1-i}\cdot e_1t^{(q)}}+\sum\limits_{i=0}^{v-1}x_2^i\otimes x_1^ux_2^{v-1-i}\cdot e_2t^{(q)}$; and
\item $I=\{2\}$: $\phi(x)=\sum\limits_{i=0}^{u-1}{x_1^ix_2^v\otimes x_1^{u-1-i}\cdot e_1e_2t^{(q)}}-[v=b-1]1\otimes x_1^u\cdot t^{(q+1)}$,\\
where $[P]=\left\{\begin{array}{lr}
1 &\text{if } P \text{ true},\\
0 &\text{if } P \text{ false}.
\end{array}\right.$
\end{itemize}
\end{Prop}
\subsection{An explicit chain map}
In the following lemmas, we will give the formula of $\tilde{f}$ based on the form of $f$ in Corollary \ref{cor.cocycles}.
\begin{Lem} \label{f_tilde1}
Let $f:F_i\to A$ be a standard cocycle of the form $(t^{(q)},x_1^ux_2^v)$ in $\mathrm{Hom}_{A^e}(F_i,A)$. For any $j \in \mathbb{N}_0$, the $A^e$-homomorphism $\tilde{f}_j:F_{i+j}\to F_j$ defined by
\begin{equation*}
\tilde{f}_j(e_Jt^{(r)})=[q \leq r]1\otimes f(t^{(q)})\cdot e_{J}t^{(r-q)}
\end{equation*}
is a lifting map of $f$ that makes \eqref{diagram1} commute.
\proof 
We shall prove this lemma by induction on $j \in \mathbb{N}_0$.
\begin{itemize}
\item $j=0$: As $f(e_Jt^{(r)})=0$ for all $e_Jt^{(r)} \neq t^{(q)}$, we then have
\begin{equation*}
\tilde{f_0}(e_Jt^{(r)})=\left\{\begin{array}{ll}
1\otimes f(t^{(q)}) &\text{if } e_Jt^{(r)} = t^{(q)},\\
0 &\text{otherwise}.
\end{array}\right.
\end{equation*}
\item $j=1$: $e_1t^{(q)}$ and $e_2t^{(q)}$ are all the basis elements of $F_{i+1}$.
\begin{align*}
\tilde{f}_1(e_1t^{(q)})&=\phi \circ \tilde{f}_0\circ d_{i+1}(e_1t^{(q)})\\
&=\phi \circ \tilde{f}_0\left((x_1\otimes 1-1\otimes x_1)\cdot t^{(q)}-e_1t^{(q-1)}\cdot d(t)\right)\\
&= 1\otimes f(t^{(q)})\cdot \phi(x_1\otimes 1-1\otimes x_1)=1\otimes f(t^{(q)})\cdot e_1.
\end{align*}
Similarly, we get $\tilde{f}_1(e_2t^{(q)})=1\otimes f(t^{(q)})\cdot e_2$.
\item Suppose that the formula holds up to   $j-1 \geq 0$. We need to show that the formula is true at $j$. Let $x=e_Jt^{(r)}$ be a basis element of degree $i+j$.\\
If $J=\emptyset$, then $r>q$. Hence, by Proposition \ref{contracting} one gets:  
\begin{align*}
\tilde{f}_j(x)&=(\phi \circ \tilde{f}_{j-1})\left(\sum\limits_{i=0}^{a-1}{x_1^i\otimes x_1^{a-1-i}\cdot e_1t^{(r-1)}}-\sum\limits_{i=0}^{b-1}{x_2^i\otimes x_2^{b-1-i}\cdot e_2t^{(r-1)}}\right)\\
&=\phi \left(\sum\limits_{i=0}^{a-1}{x_1^i\otimes x_1^{a-1-i}\cdot e_1t^{(r-1-q)}}-\sum\limits_{i=0}^{b-1}{x_2^i\otimes x_2^{b-1-i}\cdot e_2t^{(r-1-q)}}\right)\cdot 1\otimes f(t^{(q)})\\
&=1\otimes f(t^{(q)})\cdot \phi(-x_2^{b-1}\otimes 1 \cdot e_2t^{(r-1-q)})=1\otimes f(t^{(q)})\cdot t^{(r-q)}.
\end{align*}
If $J=\{1\}$, we have
\begin{align*}
\tilde{f}_j(x)&=\phi \circ \tilde{f}_{j-1}\circ d_{i+1}(e_1t^{(r)})=\phi \circ \tilde{f}_{j-1}\left(d(e_1)t^{(r)}+\sum\limits_{i=0}^{b-1}{x_2^i\otimes x_2^{b-1-i}\cdot e_1e_2t^{(r-1)}}\right)\\
&=1\otimes f(t^{(q)})\cdot \phi \left((x_1\otimes 1-1\otimes x_1)\cdot t^{(r-q)}+\sum\limits_{i=0}^{b-1}{x_2^i\otimes x_2^{b-1-i}\cdot e_1e_2t^{(r-1-q)}}\right)\\
&=1\otimes f(t^{(q)})\cdot \phi (x_1\otimes 1\cdot t^{(r-q)})\\
&=1\otimes f(t^{(q)})\cdot e_1t^{(r-q)}.
\end{align*}
By similar arguments, we will get
\begin{equation*}
\tilde{f}_j(e_2t^{(r)})=1\otimes f(t^{(q)})\cdot e_2t^{(r-q)};
\end{equation*}
and
\begin{equation*}
\tilde{f}_j(e_1e_2t^{(r)})=1\otimes f(t^{(q)})\cdot e_1e_2t^{(r-q)}.
\end{equation*}
\end{itemize}
\qed
\end{Lem} 
\begin{Lem} \label{f_tilde2}
Let $f:F_i\to A$ be a standard cocycle of the form $b(e_1t^{(q)},x_1^{u_1}x_2^{u_2})+a(e_2t^{(q)},x_1^{v_1}x_2^{v_2}))$. For $j \in \mathbb{N}_0$, the $A^e$-homomorphism $\tilde{f}_j:F_{i+j}\to F_j$ given as follows makes the diagram \eqref{diagram1} commute:
\begin{align*}
\tilde{f}_0(e_1t^{(q)})&=1\otimes bx_1^{u_1}x_2^{u_2}; \tilde{f}_0(e_2t^{(q)})=1\otimes ax_1^{v_1}x_2^{v_2};\\
\tilde{f}_{1}(t^{(q+1)})&=1\otimes bx_1^{u_1}x_2^{u_2}\cdot \delta_1e_1-1\otimes ax_1^{v_1}x_2^{v_2}\cdot \delta_2e_2;\\
\tilde{f}_{1}(e_1e_2t^{(q)})&=1\otimes ax_1^{v_1}x_2^{v_2}\cdot e_1-1\otimes bx_1^{u_1}x_2^{u_2}\cdot e_2;\\
\tilde{f}_{2j}(e_1t^{(q+j)})&=1\otimes bx_1^{u_1}x_2^{u_2}\cdot t^{(j)}-1\otimes ax_1^{v_1}x_2^{v_2}\delta_2\cdot e_1e_2t^{(j-1)};\\
\tilde{f}_{2j}(e_2t^{(q+j)})&=1\otimes ax_1^{v_1}x_2^{v_2}\cdot t^{(j)}-1\otimes bx_1^{u_1}x_2^{u_2}\cdot \delta_1e_1e_2t^{(j-1)};\\
\tilde{f}_{2j+1}(t^{(q+j+1)})&=1\otimes bx_1^{u_1}x_2^{u_2}\cdot \delta_1e_1t^{(j)}-1\otimes ax_1^{v_1}x_2^{v_2}\cdot \delta_2e_2t^{(j)};\\
\tilde{f}_{2j+1}(e_1e_2t^{(q+j)})&=1\otimes ax_1^{v_1}x_2^{v_2}\cdot e_1t^{(j)}-1\otimes bx_1^{u_1}x_2^{u_2}\cdot e_2t^{(j)},
\end{align*}
where $\delta_1=\sum\limits_{i=0}^{a-2}(i+1)x_1^{a-2-i}\otimes x_1^i$ and $\delta_2=\sum\limits_{i=0}^{b-2}(i+1)x_2^{b-2-i}\otimes x_2^i.$
%=x_2^{b-2}\otimes 1+2x_2^{b-3}\otimes x_1+\cdots +(b-1)1\otimes x_2^{b-2}
\proof
The basis of $F_i$ consists of $e_1t^{(q)}$ and $e_2t^{(q)}$. We can see that $\tilde{f}_0(e_1t^{(q)})=1\otimes bx_1^{u_1}x_2^{u_2}$ and  $\tilde{f}_0(e_2t^{(q)})=1\otimes ax_1^{v_1}x_2^{v_2}.$
In the next degree, the basis of $F_{i+1}$ consists of $t^{(q+1)}$ and $e_1e_2t^{(q+1)}$. By Proposition \ref{contracting} and the definition of $\tilde{f}$, we have that
\begin{align*}
\tilde{f}_1(t^{(q+1)})=&\phi \circ \tilde{f}_0 \circ d(t^{(q+1)})=(\phi\circ \tilde{f}_{0})\left(\sum\limits_{i=0}^{a-1}{x_1^i\otimes x_1^{a-1-i}\cdot e_1t^{(q)}}-\sum\limits_{i=0}^{b-1}{x_2^i\otimes x_2^{b-1-i}\cdot  e_2t^{(q)}}\right)\\
=&1\otimes bx_1^{u_1}x_2^{u_2}\cdot \phi\left(\sum\limits_{i=0}^{a-1}{x_1^i\otimes x_1^{a-1-i}}\right)-1\otimes ax_1^{v_1}x_2^{v_2} \cdot  \phi\left(\sum\limits_{i=0}^{b-1}{x_2^i\otimes x_2^{b-1-i}}\right)\\
=&1\otimes bx_1^{u_1}x_2^{u_2}\cdot \delta_1e_1-1\otimes ax_1^{v_1}x_2^{v_2}\cdot \delta_2e_2
\end{align*}
since
\begin{align*}
\phi\left(\sum\limits_{i=0}^{a-1}{x_1^i\otimes x_1^{a-1-i}}\right)=&\phi_1\left(1\otimes x_1^{a-1}+x_1\otimes x_1^{a-2}+\cdots +x_1^{a-1}\otimes 1\right)\\
=&\phi\left(1\otimes x_1^{a-1}\right)+\phi\left(x_1\otimes x_1^{a-2}\right)+\cdots +\phi\left(x_1^{a-1}\otimes 1\right)\\
=&\left(0\right)+\left(1\otimes x_1^{a-2}\cdot e_1\right)+\left(x_1^2\otimes x_1^{a-3}\cdot e_1+1\otimes x_1^{a-2}\cdot e_1\right)+\\
&+\cdots +\left(x_1^{a-2}\otimes 1\cdot e_1+x_1^{a-3}\otimes x_1\cdot e_1+ \cdots + 1\otimes x_1^{a-2}\cdot e_1\right)\\
=& \underbrace{\left(x_1^{a-2}\otimes 1+2x_1^{a-3}\otimes x_1+\cdots +(a-1)1\otimes x_1^{a-2}\right)}_{=:\delta_1} e_1
\end{align*}
and similarly,  $\phi\left(\sum\limits_{i=0}^{b-1}{x_2^i\otimes x_2^{b-1-i}}\right)=\delta_2e_2$.\\
Let us now consider the remaining basis element, $e_1e_2t^{(q)}$:
\begin{align*}
\tilde{f}_1(e_1e_2t^{(q)})=&\phi \circ \tilde{f}_0 \circ d(e_1e_2t^{(q)})=\phi_1 \circ \tilde{f}_0 \left(d(e_1)e_2t^{(q)}-e_1d(e_2)t^{(q)}\right)\\
=&\phi\left(\left(x_1\otimes 1-1\otimes x_1\right)1\otimes ax_1^{v_1}x_2^{v_2}-\left(x_2\otimes 1-1\otimes x_2\right)1\otimes bx_1^{u_1}x_2^{u_2}\right)\\
=&\phi\left(x_1\otimes ax_1^{v_1}x_2^{v_2}\right)-\phi_1\left( x_2\otimes bx_1^{u_1}x_2^{u_2}\right)\\
=&1\otimes ax_1^{v_1}x_2^{v_2}\cdot e_1-1\otimes bx_1^{u_1}x_2^{u_2}\cdot e_2.
\end{align*}
For the remaining degrees, we shall use induction on even and odd degrees. Suppose that the formula holds up to $2j$, we need to show the formula holds for $2j+1$. Indeed,
\begin{align*}
\tilde{f}_{2j+1}(t^{(q+j+1)})=&(\phi \circ \tilde{f}_{2j})\left(\sum\limits_{i=0}^{a-1}{x_1^i\otimes x_1^{a-1-i}\cdot  e_1t^{(q+j)}}-\sum\limits_{i=0}^{b-1}{x_2^i\otimes x_2^{b-1-i} \cdot e_2t^{(q+j)}}\right)\\
=&\phi \left( \left(\sum\limits_{i=0}^{a-1}{x_1^i\otimes x_1^{a-1-i}}\right)\cdot \left(1\otimes bx_1^{u_1}x_2^{u_2}\cdot t^{(j)}-1\otimes ax_1^{v_1}x_2^{v_2}\cdot \delta_2e_1e_2t^{(j-1)}\right)\right)\\
&-\phi \left( \left(\sum\limits_{i=0}^{b-1}{x_2^i\otimes x_2^{b-1-i}}\right)\cdot \left(1\otimes ax_1^{v_1}x_2^{v_2}\cdot t^{(j)}-1\otimes bx_1^{u_1}x_2^{u_2}\cdot \delta_1e_1e_2t^{(j-1)}\right)\right)\\
=&\phi\left(\sum\limits_{i=0}^{a-1}{x_1^i\otimes bx_1^{a-1-i} x_1^{u_1}x_2^{u_2}\cdot t^{(j)}}\right)- \phi\left(\sum\limits_{i=0}^{b-1}{x_2^i\otimes ax_2^{b-1-i}x_1^{v_1}x_2^{v_2}\cdot t^{(j)}}\right)\\
=&1\otimes bx_1^{u_1}x_2^{u_2}\cdot \delta_1e_1t^{(j)}-1\otimes ax_1^{v_1}x_2^{v_2}\cdot \delta_2e_2t^{(j)}.
\end{align*}
Using a similar argument, we get the formula of $\tilde{f}_{2j+1}(e_1e_2t^{(q+j)})$.\\
Now suppose that the formula holds up to $2j-1$, we prove that the formula at $2j$ holds.
\begin{align*}
\tilde{f}_{2j}(e_1t^{(q+j)})=&(\phi \circ \tilde{f}_{2j-1})\left(d(e_1)t^{(q+j)}+\sum\limits_{i=0}^{b-1}{x_2^i\otimes x_2^{b-1-i}}\cdot  e_1e_2t^{(q+j-1)}\right)\\
=&\phi \left(\left(x_1\otimes 1-1\otimes x_1\right)\cdot \left(1\otimes bx_1^{u_1}x_2^{u_2}\cdot \delta_1e_1t^{(j-1)}-1\otimes ax_1^{v_1}x_2^{v_2}\cdot \delta_2e_2t^{(j-1)}\right)\right)\\
&+\phi \left(\left(\sum\limits_{i=0}^{b-1}{x_2^i\otimes x_2^{b-1-i}}\right)\cdot \left(1\otimes ax_1^{v_1}x_2^{v_2}\cdot e_1t^{(j-1)}-1\otimes bx_1^{u_1}x_2^{u_2}\cdot e_2t^{(j-1)}\right)\right)\\
=&-\phi \left(x_1\otimes ax_1^{v_1}x_2^{v_2}\cdot \delta_2e_2t^{(j-1)}\right)-\phi \left(x_2^{b-1}\otimes bx_1^{u_1}x_2^{u_2}\cdot e_2t^{(j-1)}\right)\\
=&-1\otimes ax_1^{v_1}x_2^{v_2}\cdot \delta_2e_1e_2t^{(j-1)}+1\otimes bx_1^{u_1}x_2^{u_2}\cdot t^{(j)}. 
\end{align*}
Similarly we get the formula $\tilde{f}_{2j}(e_2t^{(q+j)})$ as desired.
\qed
\end{Lem}
\subsection{The cup product}
From the formula of $\tilde{f}$, the cup product can be interpreted in terms of the Yoneda product (see \cite{Wit17} Chapter 1 for more details) on $\mathrm{HH}^*(A)$ as follows. Let $f$ and $g$ be cocycles in $\mathrm{Hom}(F_i,A)$ and $\mathrm{Hom}(F_j,A)$ respectively. Then the product of these cocycles, denoted by $f\smile g$, is given by
\begin{equation*}
f\smile g:=g\circ \tilde{f}_j,
\end{equation*}
which is again a cocycle of homological degree $i+j$. Since $\tilde{f}$ is unique up to homotopy, the cup product induces a well-defined product by passing to cohomology, i.e., a multiplication on $\mathrm{HH}^*(A)$. By Lemma  \ref{f_tilde1} and \ref{f_tilde2}, we have the product of two standard residue classes in the consequence below.
\begin{Cor} \label{cup_product} 
The formula of the cup product between two standard residue classes in $\mathrm{HH}^*(A)$ is calculated as follows:
\begin{align*}
[(t^{(p)},s^\alpha)]\smile [(t^{(q)},s^\beta)]&=[(t^{(p+q)},s^{\alpha +\beta})];\\
[(t^{(p)},s^\alpha)]\smile \left[b(e_1t^{(q)},s^{\beta-m_1})\right.+\left.a(e_2t^{(q)},s^{\beta-m_2})\right]&=[b(e_1t^{(p+q)},s^{\alpha+\beta-m_1})+a(e_2t^{(p+q)},s^{\alpha+\beta-m_2})];\\
\left[b(e_1t^{(p)},s^{\alpha-m_1})+a(e_2t^{(p)},s^{\alpha-m_2})\right]&\smile 
\left[b(e_1t^{(q)},s^{\beta-m_1})+a(e_2t^{(q)},s^{\beta-m_2})\right]=0.
\end{align*}
Moreover, the multiplication is commutative.
\proof
Let $f:=(t^{(p)},s^\alpha)$ and $g:=(t^{(q)},s^\beta)$. We calculate the first product as follows:
\begin{align*}
(f\smile g)(e_Jt^{(u)})=&g\circ \tilde{f}(e_Jt^{(u)})=[p\leq u]g(e_Jt^{(u-p)})\cdot 1\otimes f(t^{(p)})\\
=&\left\{\begin{array}{ll}
{g(t^{(q)})\cdot f(t^{(p)})}&{\text{if }e_Jt^{(u-p)}=t^{(q)}}\\
{0}&{\text{otherwise}}
\end{array}\right.
=\left\{\begin{array}{ll}
{s^{\alpha+\beta}}&{\text{if }e_Jt^{(u)}=t^{(p+q)},}\\
{0}&{\text{otherwise}.}
\end{array}\right.
\end{align*}
For the second one, let $f:=(t^{(p)},s^\alpha)$ and $g:=b(e_1t^{(q)},s^{\beta-m_1})+a(e_2t^{(q)},s^{\beta-m_2})$. By a similar computation, we get the second product:
\begin{align*}
(f\smile g)(e_Jt^{(u)})=&\left\{\begin{array}{ll}
{g(e_1t^{(q)})\cdot f(t^{(p)})}&{\text{if }e_Jt^{(u-p)}=e_1t^{(q)},}\\
{g(e_2t^{(q)})\cdot f(t^{(p)})}&{\text{if }e_Jt^{(u-p)}=e_2t^{(q)},}\\
{0}&{\text{otherwise},}
\end{array}\right.\\
=&\left\{\begin{array}{ll}
{bs^{\alpha+\beta-m_1}}&{\text{if }e_Jt^{(u)}=e_1t^{(p+q)},}\\
{as^{\alpha+\beta-m_2}}&{\text{if }e_Jt^{(u)}=e_2t^{(p+q)},}\\
{0}&{\text{otherwise}.}
\end{array}\right.
\end{align*}
Now we take $f=b(e_1t^{(p)},s^{\alpha-m_1})+a(e_2t^{(p)},s^{\alpha-m_2})$ of degree $i$ and $g=b(e_1t^{(q)},s^{\beta-m_1})+a(e_2t^{(q)},s^{\beta-m_2})$ of degree $j$. The basis of $F_{i+j}$ consists of $e_1e_2t^{(p+q)}$ and $t^{(p+q+1)}$. Apply Lemma \ref{f_tilde2}, replace $x_1, x_2$ by $s^b, s^a$ and notice that $j=2q+1$, we get that:
\begin{align*}
(f\smile g)(e_1e_2t^{(p+q)})=&g\left(\tilde{f}_j(e_1e_2t^{(p+q)})\right)=g(1\otimes as^{\alpha-m_2}\cdot e_1t^{(q)}-1\otimes bs^{\alpha-m_1}\cdot e_2t^{(q)})\\
=&(1\otimes as^{\alpha-m_2})\cdot bs^{\beta-m_1}-(1\otimes bs^{\alpha-m_1})\cdot as^{\beta-m_2}=0.
\end{align*}
Let us consider the remaining basis element:
\begin{align*}
(f\smile g)(t^{(p+q+1)})=&g(1\otimes bs^{\alpha-m_1}\delta_1e_1t^{(q)}-1\otimes as^{\alpha-m_2}\delta_2e_2t^{(q)})\\
=&\dfrac{a(a-1)}{2}b^2s^{\alpha+\beta-2m_1+b(a-2)}-\dfrac{b(b-1)}{2}a^2s^{\alpha+\beta-2m_2+a(b-2)}\\
=&\dfrac{ab(a-b)}{2}s^{\alpha+\beta-ab}.
\end{align*}
Here, we consider $\delta_1$, $\delta_2$ using variable $s$. So we have showed that $$f\smile g=\dfrac{ab(a-b)}{2}(t^{(p+q+1)},s^{\alpha+\beta-ab}).$$
Now we will state that $(t^{(p+q+1)},s^{\alpha+\beta-ab})$ belongs to the image of $\partial$, i.e., its residue class in $\mathrm{HH}^*(A)$ is zero. By Remark \ref{non_zero_cocycles}, we will show that $\alpha+\beta-ab-m_1\in S$ or $\alpha+\beta-ab-m_2\in S$. From Corollary \ref{cor.cocycles}, there are two options for $\alpha-m_1$ and  $\beta-m_1$, which are $m_2$ and $\gamma+b$ for some $\gamma \in S$.
\begin{itemize}
\item If $\alpha-m_1=\beta-m_1=m_2$, then $\alpha=\beta=m_1+m_2$. Hence, $\alpha+\beta-ab-m_1=m_1+a(b-2)\in S$ and  $\alpha+\beta-ab-m_2=m_2+(a-2)b\in S$. 
\item If $\alpha-m_1=m_2$ and $\beta-m_1=\gamma+b$ for some $\gamma \in S$, then $\alpha+\beta-ab-m_1=\gamma+m_2\in S$ and $\alpha+\beta-ab-m_2=\gamma+m_1\in S$. Similarly for $\alpha-m_1=\gamma+b$ and $\beta-m_1=m_2$.
\item If $\alpha-m_1=\gamma+b$ and $\beta-m_1=\eta+b$ for some $\gamma, \eta \in S$, then $\alpha+\beta-ab-m_1=\gamma +\eta +b\in S$ and $\alpha+\beta-ab-m_2=\gamma +\eta +a\in S$.
\end{itemize}
\qed
\end{Cor}
By supplementing the module $\mathrm{HH}^*(A)$ with a multiplicative structure, this module becomes a $k$-algebra. 
By Corollary \ref{cor.cocycles} and \ref{cup_product}, we have the description of the generators for the algebra $\mathrm{HH}^*(A)$ as follows.
\begin{Rem} \label{generators}
(a) For the basic element of the form $[(t^{(p)},s^\alpha)]$ (where $\alpha \in S$) there are $u, v \in \mathbb{N}_0$ such that $\alpha = ua+vb$. Then we can write $[(t^{(p)},s^\alpha)]$ as a product of the elements $[(t,1)]$, $[(1,s^a)]$ and $[(1,s^b)]$.\\
(b) Likewise, a basic element of the form $[b(e_1t^{(q)},s^{\alpha-m_1})+a(e_2t^{(q)},s^{\alpha-m_2})]$ (where $\alpha \in S$) is written as a product of $[(t,1)]$, $[(1,s^a)$, $(1,s^b)]$ and either $[b(e_1,s^{m_2})+a(e_2,s^{m_1})]$ or $[b(e_1,s^{b})+a(e_2,s^{a})]$, where the two last elements occur once for such a basic element of this type.
\end{Rem}
Now we are in the position to give a proof of Theorem A which we stated in the introduction.
\subsection{Proof of Theorem A}\label{proofA}
The Hochschild cohomology module $\mathrm{HH}^*(A)$ consists of the cosets of the cocycles in $\mathrm{Ker}(\partial)$. We set $X_1$ to be the element $(1,s^a)$. Similarly, we have $X_2$ for $[(1,s^b)]$, $Y_1$ for $[b(e_1,s^{b})+a(e_2,s^{a})]$, $Y_2$ for $[b(e_1,s^{m_2})+a(e_2,s^{m_1})]$ and $T$ for $[(t,1)]$. Let us introduce a multidegree `$\mathrm{mdeg}$' combined from an $\mathbb{N}_0$-grading (on the first argument) and a $\mathbb{Z}$-weight (on the second argument) by setting $\mathrm{mdeg}(e_1,1)=(1,-b)$, $\mathrm{mdeg}(e_2,1)=(1,-a)$, $\mathrm{mdeg}(1,s)=(0,1)$, $\mathrm{mdeg}(t,1)=(2,-d_1d_2)$. Then we consider the decomposition of $\mathrm{HH}^*(A)$ induced by our multidegree. The differential $\partial$ is a 1-homogeneous morphism with respect to the grading and a 0-homogeneous morphism with respect to the weight. By Remark \ref{generators}, we know that $\mathrm{HH}^*(A)$ is generated by $X_1$, $X_2$, $Y_1$, $Y_2$ and $T$. The degree (`$\mathrm{deg}$') and the weight (`$\mathrm{wt}$') of these elements follow from the multidegree. To show that the relations in the theorem are satisfied, we use Corollary \ref{cup_product} as follows.
\begin{itemize}
\item As $(1,s^a)^b=(1,s^{ab})=(1,s^b)^a$, we have the first relation, $X_1^b-X_2^a$.
\item Using the formula in Corollary \ref{cup_product}, we have the relation $Y_1^2$,$Y_2^2$ and $Y_1Y_2$.
\item By Remark \ref{non_zero_cocycles}, the standard cocycles in the image of $\partial$ consist of:\\
$(t^{(q)},s^{\alpha})$, where $q>0$, $\alpha \in S$, $\alpha - {m_1} \in S$; \\
$(t^{(q)},s^{\alpha})$, where $q>0$, $\alpha \in S$,  $\alpha - {m_2} \in S$; and\\
$b(e_1t^{(q)},s^{\alpha +m_2})+a(e_2t^{(q)},s^{\alpha +m_1})$, where $q>0$, $\alpha \in S$, $\alpha - F(S) \in S$.\\
From this, we can deduce the relations $X_1^{b-1}T$, $X_2^{a-1}T$ and $Y_2T$.
\end{itemize}
So far, we have obtained all generators and relations displayed in the statement. Now we will prove that there is an isomorphism between the algebras, $\mathrm{HH}^*(A)$ and $k[X_1,X_2,Y_1,Y_2,T]/\mathcal{I}$, by showing that there is a bigraded bijection between a $k$-basis of each.

We first describe the $k$-basis of the algebra $k[X_1,X_2,Y_1,Y_2,T]/\mathcal{I}$. The Gr\"obner basis of $\mathcal{I}$ with respect to the pure lexicographic term order $X_1 \prec X_2 \prec Y_1 \prec Y_2 \prec T$ on $k[X_1,X_2,Y_1,Y_2,T]$ is determined as follows:
$X_2^a-X_1^b$, $X_1^{b-1}T$, $X_2^{a-1}T$, $Y_2T$, $Y_1^2$, $Y_2^2$, $Y_1Y_2$, $X_1Y_2-X_2^{a-1}Y_1$, $X_2Y_2-X_1^{b-1}Y_1$. The leading terms of this  Gr\"obner base are $X_2^a$, $X_1^{b-1}T$, $X_2^{a-1}T$, $Y_2T$, $Y_1^2$, $Y_2^2$, $Y_1Y_2$, $X_1Y_2$, $X_2Y_2$. From here, one has a $k$-basis of the algebra $k[X_1,X_2,Y_1,Y_2,T]/\mathcal{I}$ consisting of the following elements:
\begin{itemize}
\item $X_1^uX_2^v$, where $u \geq 0$, $0 \leq v < a$;
\item $X_1^uX_2^vT^q$, where $0 \leq u < b-1$, $0 \leq v < a-1$, $q>0$;
\item $X_1^uX_2^vY_1$, where $u \geq 0$, $0 \leq v < a$;
\item $X_1^uX_2^vY_1T^q$, where $0 \leq u < b-1$, $0 \leq v < a-1$, $q>0$; and
\item $Y_2$.
\end{itemize}

The $k$-basis of  $\mathrm{HH}^*(A)$ was described in Remark \ref{non_zero_cocycles}. It can be easily seen that there is a bigraded one-to-one correspondence between: $X_1^uX_2^v$, where $u \geq 0$, $0 \leq v < a$ and $(1, s^\alpha)$, where $\alpha \in S$; $Y_2$ and $b(e_1,s^{m_2})+a(e_2,s^{m_1})$; $X_1^uX_2^vY_1$, where $u \geq 0$, $0 \leq v < a$ and $b(e_1,s^{\alpha +b})+a(e_2,s^{\alpha +a})$, where $\alpha \in S$. Now we will show that the rest of the $k$-bases of  $\mathrm{HH}^*(A)$ and $k[X_1,X_2,Y_1,Y_2,T]/\mathcal{I}$ are corresponding to each other as well.

(i) $X_1^uX_2^vT^q$, where $0 \leq u < b-1$, $0 \leq v < a-1$, $q>0$ and $(t^{(q)},s^\alpha)$, where $q>0$, $\alpha \in S$, $\alpha -m_1 \notin S$, $\alpha - m_2 \notin S$ are equivalent. Indeed, suppose that $\alpha = ua+vb$ ($u,v \in \mathbb{N}_0$). We will show that
$$
\left\{\begin{array}{l}
{u<b-1}\\
{v<a-1}
\end{array}
\right.
\Leftrightarrow 
\left\{\begin{array}{l}
{ua+vb-m_1 \notin S}\\
{ua+vb-m_2 \notin S}
\end{array}
\right..
$$
``$\Leftarrow$" Suppose on the contrary that $u \geq b-1$. Then,  $ua+vb-m_2=ua+vb-a(b-1)=vb+a\left(u-(b-1)\right)\in S$, which is a contradiction. Similar for $v$. ``$\Rightarrow$" By Lemma \ref{p2lemma1}, $v<a-1$  implies that $vb-a \notin S$. Since $u<b-1$, $\gamma=ua+vb-m_2=ua+vb-a(b-1)=vb+a\left(u-(b-1)\right)\notin S$. If not, $\gamma \in S$, so $vb-a=\gamma+a(b-2-u)\in S$, which is inconsequential.

(ii) $X_1^uX_2^vY_1T^q$, where $0 \leq u < b-1$, $0 \leq v < a-1$, $q>0$ corresponds to $b(e_1t^{(q)},s^{\alpha +b})+a(e_2t^{(q)},s^{\alpha +a})$, where $\alpha \in S$, $\alpha-F(S) \notin S$, $q>0$. Suppose that $\alpha = ua+vb$ ($u,v \in \mathbb{Z}, \geq 0$). We will show that
$$
\left\{\begin{array}{l}
{u<b-1}\\
{v<a-1}
\end{array}
\right.
\Leftrightarrow 
ua+bv-F(S)\notin S.
$$
``$\Leftarrow$" Suppose on the contrary that $u \geq b-1$ or $v \geq a-1$. Then we have $ua+bv-ab+a+b=(u-b+1)a+bv+b \in S$ or $ua+bv-ab+a+b=(v-a+1)b+au+a \in S$, which contradicts $ua+bv-F(S)\notin S$. ``$\Rightarrow$" Suppose that $u<b-1$ and $v<a-1$. We need to show that $ua+bv-F(S)\notin S$. If $ua+bv-F(S)\in S$, then $ua+bv-F(S)=-a+(v+1)b+(u-b+2)a \in S$, where $u-b+2\leq 0$. This implies that  $-a+(v+1)b \in S$ (where $v+1<a$) which is impossible by Lemma \ref{p2lemma1}.

Hence, we have proved that the Hochschild cohomology ring $\mathrm{HH}^*(A)$ is isomorphic to the quotient ring $k[X_1,X_2,Y_1,Y_2,T]/\mathcal{I}$. \qed
\subsection{The Hilbert series of $\mathrm{HH}^*(A)$} \label{hilbert1}
Let $\mathbf{H}_{m,n}$ be the $k$-module generated by the elements whose degree is $(m,n)\in \mathbb{N}_0\times \mathbb{Z}$. The Hilbert series of $\mathrm{HH}^*(A)=\mathop \oplus\limits_{(m,n)\in \mathbb{N}_0\times \mathbb{Z}}\mathbf{H}_{m,n}$ as an $\mathbb{N}_0\times \mathbb{Z}$-graded vector space via the grading above is the formal series:
$$\mathcal{H}_{\mathrm{HH}^*(A)}(x,y)=\mathop\sum\limits_{(m,n)\in \mathbb{N}_0\times \mathbb{Z}}\mathrm{dim}_k(\mathbf{H}_{m,n})x^my^n.$$
This series is computed based on the Hilbert series of the non-zero cocycles, whose description is listed in Remark \ref{non_zero_cocycles}.  We will use the decomposition introduced in Theorem \ref{main.theorem} in computing the Hilbert series.

(i) We have that $(1,s^{\alpha})$, where $\alpha \in S$, contributes the series
$$\mathrm{H}_1=\mathcal{H}_{k[s^a,s^b]}(x,y)=\dfrac{1-y^{ab}}{(1-y^{a})(1-y^{b})}.$$

(ii) Now we consider the non-zero cocycles of type \eqref{cocycles_type2} in Remark \ref{non_zero_cocycles}, $(t^{(q)},s^{\alpha})$, where $q>0$, $\alpha \in S$, $\alpha -m_1\notin S$ and $\alpha -m_2\notin S$. The element $(t^{(q)},s^{\alpha})$,where $\alpha \in S$, has degree $\left(2q,q\left(-ab\right) +\alpha\right)$. This element contributes the term $x^{2q}y^{q\left(-ab\right) +\alpha}$, which is $(x^{2}y^{-ab})^qy^{\alpha}$ equivalently. We notice that
\begin{equation*}
\{\alpha \in S \ | \ \alpha-m_1 \notin S \text{ and } \alpha-m_2  \notin S\} = S \smallsetminus (S_{1} \cup S_{2})
\end{equation*}
and, by the principle of inclusion-exclusion, we have that
\begin{equation} \label{PIE}
\sum \limits_{\alpha \in S \smallsetminus (S_{1} \cup S_{2})}{y^{\alpha}}=\sum \limits_{\alpha \in S}{y^{\alpha}}-\left(\sum \limits_{\alpha \in S_{1}}{y^{\alpha}}+\sum \limits_{\alpha \in S_2}{y^{\alpha}}\right)+\sum \limits_{\alpha \in S_1\cap S_2}{y^{\alpha}}. 
\end{equation}
By Lemma \ref{lemma2}, we already have the detailed description of $S_1$, $S_2$ and $S_{1}\cap S_{2}$. Now we are able to calculate the Hilbert series formed by this kind of cohomology classes.
\begin{itemize}
\item  The series given by all $(t^{(q)},s^{\alpha})$, where $\alpha \in S$ and $q>0$ is
$$ \mathrm{H}_{21}=\dfrac{x^2y^{-ab}}{1-x^2y^{-ab}}\cdot\mathrm{H}_1.$$
\item The element $(t^{(q)},s^{\alpha})$, where $\alpha \in S_1$ is written as $(t^{(q)},s^{\gamma+m_1})$, where $\gamma \in S$. Hence, the corresponding degree is $\left(2q,q\left(-ab\right) + m_1 +\gamma \right)$, which contributes the term $x^{2q}y^{q\left(-ab\right) + m_1 +\gamma}$. Then, the series given by all $(t^{(q)},s^{\gamma+m_1})$, where $\gamma \in S$, $q>0$, is
$$ \mathrm{H}_{22}=\dfrac{x^2y^{-ab}}{1-x^2y^{-ab}}\cdot y^{m_1}\cdot \mathrm{H}_1.$$
\item Similarly, the series given by all $(t^{(q)},s^{\alpha})$, where $\alpha \in S_2$, $q>0$, is
$$ \mathrm{H}_{23}=\dfrac{x^2y^{-ab}}{1-x^2y^{-ab}}\cdot y^{m_2}\cdot \mathrm{H}_1.$$
\item The element $(t^{(q)},s^{\alpha})$, where $\alpha \in S_1\cap S_2$, is $(t^{(q)},s^{m_1+m_2})$ or $(t^{(q)},s^{ab+\gamma})$, where $\gamma \in S$. Hence, by a similar argument, we find out that the series for these elements is
$$ \mathrm{H}_{24}=\dfrac{x^2y^{-ab}}{1-x^2y^{-ab}}\cdot \left(y^{m_1+m_2}+y^{ab}\cdot \mathrm{H}_1\right).$$
\end{itemize}
By \eqref{PIE}, the Hilbert series for the elements of type \eqref{cocycles_type2} is $$\mathrm{H}_2=\mathrm{H}_{21}-(\mathrm{H}_{22}+\mathrm{H}_{23})+\mathrm{H}_{24}.$$

(iii) For the cocycles of type \eqref{cocycles_3} in Remark \ref{non_zero_cocycles},  we have the single cocycle $b(e_1,s^{m_2})+a(e_2,s^{m_1})$ and the cocycles $b(e_1t^{(q)},s^{b+\alpha})+a(e_2t^{(q)},s^{a+\alpha})$, where $\alpha \in S$ and if $q>0$, $\alpha-F(S) \notin S$.
\begin{itemize}
\item The element $b(e_1,s^{m_2})+a(e_2,s^{m_1})$ of degree $(1,ab-a-b)$ and the elements $b(e_1,s^{b+\alpha})+a(e_2,s^{a+\alpha})$ of degree $(1,\alpha)$ contribute the series
$$\mathrm{H}_{31}=xy^{ab-a-b}+x\cdot \mathrm{H}_1.$$
\item For the remaining elements, notice that
\begin{equation*}
\{\alpha\in S \ | \ \alpha-F(S) \notin S\}=S \smallsetminus \{\gamma + F(S) \ | \ \gamma \in S \smallsetminus \{0\} \}.
\end{equation*} 
So we have the series
$$\mathrm{H}_{32}=x\cdot \dfrac{x^2y^{-ab}}{1-x^2y^{-ab}} \cdot\mathrm{H}_1,$$
which corresponds to the $b(e_1t^{(q)},s^{b+\alpha})+a(e_2t^{(q)},s^{a+\alpha})$, where $q>0$ and $\alpha \in S$.\\
And the series
$$\mathrm{H}_{33}=xy^{ab-a-b}\cdot \dfrac{x^2y^{-ab}}{1-x^2y^{-ab}} \cdot(\mathrm{H}_1-1),$$
which corresponds to the elements $b(e_1t^{(q)},s^{b+\gamma+F(S)})+a(e_2t^{(q)},s^{a+\gamma+F(S)})$, where $q>0$ and $\gamma \in S\smallsetminus \{0\}$.\\
\end{itemize}
Now we get the Hilbert series for all elements of type \eqref{cocycles_3}, which is $\mathrm{H}_3=\mathrm{H}_{31}+\mathrm{H}_{32}-\mathrm{H}_{33}$. Hence, the Hilbert series for Case I is: $\mathcal{H}_{\mathrm{HH}^*(A)}(x,y)=\mathrm{H}_{1}+\mathrm{H}_{2}+\mathrm{H}_{3}$.
\section{The ring structure of $\mathrm{HH}^*(A)$ - Case II}\label{s4}
In this section, we will use the same arguments as in Case I to describe the ring structure of  $\mathrm{HH}^*(A)$ in the case that $\mathrm{char}(k)$ is a divisor of $a$. Some of our results shall be stated without proof because the reader can establish them analogously to the previous case.
\subsection{Formula of the cup product}
Since $\mathrm{char}(k)$ is a divisor of $a$, the formula of $\partial$ becomes
\begin{equation*} \partial(e_It^{(q)},s^{\alpha})=
\left\{
	\begin{array}{ll}
		0  & \text{if }I = \emptyset \text{ or }\{1\}, \\
		-b(t^{(q+1)}, s^{\alpha+m_2})  & \text{if }I = \{2\}, \\
		b(e_1t^{(q+1)}, s^{\alpha+m_2})  & \text{if }I = \{1,2\}.
	\end{array}
\right. 
\end{equation*}
Then we have an immediate consequence of our information on the kernel and the image of $\partial$ as follows.
\begin{Cor} \label{cocycles_2}
\begin{itemize}
\item [(i)] The kernel of $\partial$ is spanned by $(e_It^{(q)},s^\alpha)$, where $\alpha \in S$ and $I = \emptyset$ or $I=\{1\}$.
\item [(ii)] The image of $\partial$ is spanned by $(e_It^{(q)},s^\alpha)$, where $\alpha \in S$, $I = \emptyset$ or $\{1\}$, $\alpha-m_2\in S$ and $q>0$.
\end{itemize}
\end{Cor}
As the explicit chain map was constructed independently from the characteristic $n$, we can interpret this chain map from Case I for Case II.
\begin{Lem}
\begin{itemize}
\item[(i)] Let $f:F_i\to A$ be a cocycle of the form $(t^{(q)},s^\alpha)$ in $\mathrm{Hom}_{A^e}(F_i,A)$. For any $j \in \mathbb{N}_0$, the $A^e$-homomorphism $\tilde{f}_j:F_{i+j} \to F_j$ given by $$\tilde{f}_j(e_Jt^{(r)})=[q\leq r]e_Jt^{(r-q)}\cdot1 \otimes f(t^{(q)})$$
is a lifting of $f$ that makes \eqref{diagram1} commute.
\end{itemize}
\item[(ii)] If $f:F_i\to A$ is a cocycle of the form $(e_1t^{(q)},s^\alpha)$ in $\mathrm{Hom}_{A^e}(F_i,A)$, then for any $j \in \mathbb{N}_0$, the $A^e$-homomorphism $\tilde{f}_j:F_{i+j} \to F_j$ is given by:
\begin{equation*}
\tilde{f}_0(e_1t^{(q)})=1\otimes s^\alpha; \tilde{f}_0(e_2t^{(q)})=0;
\tilde{f}_{1}(t^{(q+1)})=1\otimes s^\alpha\delta_1e_1; \tilde{f}_{1}(e_1e_2t^{(q)})=-1\otimes s^\alpha e_2;
\end{equation*}
\begin{equation*}
\tilde{f}_{2j}(e_1t^{(q+j)})=1\otimes s^\alpha t^{(j)};
\tilde{f}_{2j}(e_2t^{(q+j)})=-1\otimes s^\alpha \delta_1e_1e_2t^{(j-1)};
\end{equation*}
\begin{equation*}
\tilde{f}_{2j+1}(t^{(q+j+1)})=1\otimes s^\alpha\delta_1e_1t^{(j)};
\tilde{f}_{2j+1}(e_1e_2t^{(q+j)})=-1\otimes s^\alpha e_2t^{(j)}.
\end{equation*}
\end{Lem}
\begin{Cor} \label{cup_product2}
The formula of the cup product between two standard residue classes in $\mathrm{HH}^*(A)$ is given by:
\begin{equation*}
[(t^{(p)},s^\alpha)]\smile [(t^{(q)},s^\beta])=[(t^{(p+q)},s^{\alpha+\beta})]
\end{equation*}
\begin{equation*}
[(t^{(p)},s^\alpha)]\smile [(e_1t^{(q)},s^\beta)]=[(e_1t^{(p+q)},s^{\alpha+\beta})]
\end{equation*}
\begin{equation*}
[(e_1t^{(p)},s^\alpha)]\smile[(e_1t^{(q)},s^\beta)]=\left\{\begin{array}{ll}
{[(t^{(p+q+1)},s^{\alpha+\beta+b(a-2)})]}& {\text{ if } \mathrm{char}(k)=2 \text{ and } 4 \nmid a,}\\
{0} &{\text {otherwise}.}
\end{array}\right.
\end{equation*}
\proof
The two first formulas are obtained by computing directly. For the last formula, we have 
\begin{equation*}
(e_1t^{(p)},s^\alpha)\smile(e_1t^{(q)},s^\beta)=\dfrac{a(a-1)}{2}(t^{(p+q+1)},s^{\alpha+\beta+b(a-2)}).
\end{equation*}
Recall that $\mathrm{char}(k)$ is a divisor of $a$. If $\mathrm{char}(k)\neq 2$ or $\mathrm{char}(k)=2$ and $a$ is divisible by 4, then $\mathrm{char}(k)$ is a divisor of $\dfrac{a(a-1)}{2}$. Hence, we have $\dfrac{a(a-1)}{2}=0$. If $\mathrm{char}(k)=2$ and 4 is not a divisor of $a$, then $a=2n$ where $n$ is an odd number. Then we get $\dfrac{a(a-1)}{2}=n(2n-1)\equiv 1$ modulo 2.
\qed
\end{Cor}
\subsection{Proof of Theorem B}\label{proofB}
All cocycles are combinations of the elements $(t^{(q)},s^\beta)$ and $(e_1t^{(q)},s^\beta)$ where $\beta \in S$. By Corollary \ref{cup_product2}, we can see that all basis cocycles are products of $(1,s^a)$, $(1,s^b)$, $(e_1,1)$ and $(t,1)$. So we set $X_1$ to be the coset of the element $(1,s^a)$, $X_2$ for $(1,s^b)$, $Y$ for $(e_1,1)$ and $T$ for $(t,1)$. Then these are generators of the ring. In addition, we easily obtain all the relations $X_1^{b}-X_2^{a}$ (as $(1,s^a)^b=(1,s^b)^a$), $X_1^{b-1}T$ by Corollary \ref{cocycles_2} (ii) and $Y^2-X_2^{a-2}T$ if $\mathrm{char}(k)=2$ and $4 \nmid a$ (or $Y^2$ otherwise) by Corollary \ref{cup_product2}.

Now we will show that there is a bigraded bijection between the $k$-bases of $k[X_1,X_2,Y,T]/\mathcal{I}$ and $\mathrm{HH}^*(A)$. Let us start with $k[X_1,X_2,Y,T]/\mathcal{I}$.  The Gr\"obner basis of $\mathcal{I}$ with respect to the pure lexicographic term order $Y \succ X_1 \succ X_2 \succ T$ consists of $X_1^{b}-X_2^{a}$, $Y^2-X_2^{a-2}T$, $X_1^{b-1}T$ and $X_2^aT$ in the case that $\mathrm{char}(k)=2$ and $a$ is not divisible by 4. The other case is very similar. Moreover, the Gr\"obner basis has the same leading terms with respect to the above order, so we can skip this case. Then, we get the $k$-basis of $k[X_1,X_2,Y,T]/\mathcal{I}$ as follows:
\begin{itemize}
\item $X_1^uX_2^vY^i$, where  $0\leq u <b$, $v \geq 0$ and  $i\in \{0,1\}$;
\item $X_1^uX_2^vY^iT^q$, where $0\leq u <b-1$, $0 \leq v <a$, $i\in \{0,1\}$ and $q>0$.
\end{itemize}
By Corollary \ref{cocycles_2}, we can deduce that the basis cocyles corresponding to the non-zero elements in $\mathrm{HH}^*(A)$ are
$(t^{(q)},s^\alpha)$ and $(e_1t^{(q)},s^\alpha)$, where $\alpha \in S$ and if $q>0$, $\alpha-m_2\notin S$. In the following, we will see the correspondence between the $k$-bases of the two rings:
\begin{itemize}
\item $X_1^uX_2^v$ ($0\leq u <b$, $v \geq 0$) corresponds to $(1,s^\alpha)$ where $\alpha \in S$.
\item $X_1^uX_2^vY$ ($0\leq u <b$, $v \geq 0$) corresponds to $(e_1,s^\alpha)$ where $\alpha \in S$.
\item To show that $X_1^uX_2^vY^iT^q$ ($0\leq u <b-1$, $0 \leq v <a$, $i\in \{0,1\}$ and $q>0$) corresponds to $(e_It^{(q)},s^\alpha)$ ($\alpha \in S$, $\alpha-m_2\notin S$, $I=\emptyset$ or $I=\{1\}$), we have to prove that 
\begin{equation*}
ua+vb-a(b-1) \notin S \Leftrightarrow \left\{
\begin{array}{ll}
{0\leq u <b-1}\\
{0 \leq v <a}
\end{array}
\right.,
\end{equation*}
where $\alpha=ua+vb$, $u,v \geq 0$.\\
Indeed, if $u \geq b-1$ or $v\geq a$, then  $ua+vb-a(b-1)\in S$, which is a contradiction. For the other implication, suppose that we have the hypothesis on the right hand side, i.e., we can write $u=b-2-d$ and $v=a-1-e$ for some $d,e \geq 0$. If we have $ua+vb-a(b-1) \in S$, then $(b-2-d)a+(a-1-e)b-a(b-1)=F(S)-ad-be \in S$. This implies that $F(S)\in S$, which is a contradiction again.\qed
\end{itemize}

\subsection{The Hilbert series of $\mathrm{HH}^*(A)$}\label{hilbert2}
We define the formal series as for the previous case and use the same decomposition for grading the Hochschild cohomology ring $\mathrm{HH}^*(A)$ in this case.	Using a similar argument to Case I, we now get the Hilbert series for $\mathrm{HH}^*(A)$ in this case as follows:
$$\mathcal{H}_{\mathrm{HH}^*(A)}(x,y)=\left(1 + \dfrac{x^2y^{-ab}}{1-x^2y^{-ab}}\cdot (1-y^{m_2})\right)\cdot (1+ xy^{-a})\mathrm{H}_1,$$
 where $\mathrm{H}_1:=\mathcal{H}_{k[s^a,s^b]}(x,y)$.
 

\begin{thebibliography}{30}

\bibitem{Hol00}  Holm, Th. (2000). Hochschild cohomology rings of algebras $k\left[X\right]/\left(f\right)$, {\it Contributions to Algebra and Geometry}, 41(1):291$-$301.

\bibitem{TS18} Tran, N, Sk\"oldberg, E. (2018). Hochschild cohomology of square-free monomial complete intersections. 
To appear in Communications in Algebra.

\bibitem{GGRV} Guccione, J. A.,  Guccione, J. J.,  (1991). Hochschild Homology of complete intersections. {\it Journal of Pure and Applied Algebra}, 74:159$-$176.

\bibitem{ES} Sk\"oldberg, E. (2005). Morse theory from an algebraic viewpoint. \textit{Trans. AMS}, 358(1): 115$-$129.

\bibitem{Wit17} Witherspoon, S. (2017). {\it An introduction to Hochschild Cohomology}, Texas A\&M University.

\end{thebibliography}
\end{document}